\documentclass[num-refs]{wiley-article}
\usepackage{siunitx}
\usepackage{graphicx}
\usepackage{amssymb}
\usepackage{amsmath}
\usepackage{dcolumn}
\usepackage{natbib}
\newcommand{\norm}[1]{\left\lVert #1 \right\rVert}

\papertype{Original Article}
\paperfield{Process Systems Engineering}

\title{Some manifold learning considerations \\
towards explicit model predictive control}

\author[1\authfn{1}]{Robert~J.~Lovelett}
\author[1]{Felix~Dietrich}
\author[1]{Seungjoon~Lee}
\author[1]{Ioannis~G.~Kevrekidis}

\affil[1]{Department of Chemical and Biomolecular Engineering, Johns Hopkins University, Baltimore, MD 21218, USA}

\fundinginfo{This work was partially supported by DARPA (Lagrange Program, Contract No. N6601-18-C-4031, as well as the Physics of Artificial Intelligence Program, Agreement HR00111890032).}

\presentadd[\authfn{1}]{Department of Chemical and Biological Engineering, Princeton University, Princeton, NJ 08544, USA}

\corraddress{Robert J. Lovelett, Department of Chemical and Biomolecular Engineering, Johns Hopkins University, Baltimore, MD 21218, USA}
\corremail{lovelett@princeton.edu}

\corraddress{Ioannis G. Kevrekidis, Department of Chemical and Biomolecular Engineering, Johns Hopkins University, Baltimore, MD 21218, USA}
\corremail{yannisk@jhu.edu}

\runningauthor{Robert J. Lovelett et al.}
\begin{document}

\maketitle

\begin{abstract}
%% Text of abstract

Model predictive control (MPC) is a \emph{de facto} standard control algorithm across the process industries. There remain, however, applications where MPC is impractical because an optimization problem is solved at each time step. We present a link between explicit MPC formulations and manifold learning to enable facilitated prediction of the MPC policy. Our method uses a similarity measure informed by control policies and system state variables, to ``learn'' an intrinsic parametrization of the MPC controller using a diffusion maps algorithm, which will also discover a low-dimensional control law when it exists as a smooth, nonlinear combination of the state variables. We use function approximation algorithms to project points from state space to the intrinsic space, and from the intrinsic space to policy space. The approach is illustrated first by ``learning'' the intrinsic variables for MPC control of constrained linear systems, and then by designing controllers for an unstable nonlinear reactor.

\end{abstract}

%% main text
\section{Introduction}
\label{sec:intro}

Model predictive control (MPC) is a \textit{de facto} standard method for control across the process industries \cite{Qin2003ATechnology}.
In MPC, which is also known as receding horizon control, the control action is calculated on-line by solving a receding horizon optimal control problem at each time step to determine the subsequent control action to take \cite{Garcia1989,Mayne2000,Mayne2014}.
MPC for constrained linear systems has been well-established for some time now \cite{Mayne2000}, and that success has led to a significant effort to extend MPC to systems that are harder to control due to stochasticity \cite{Mayne2014}, nonlinearity \cite{Meadows1996}, decentralization \cite{Christofides2013}, or other challenges.
Many of the theoretical challenges associated with these more complex MPC problems have been overcome \cite{Meadows1996}, enabling practitioners to use MPC in new applications that are even more demanding.
Because, however, MPC entails solving an optimization problem at every sampled time step, it will always be limited by the computational time it takes to solve that problem.
Due to limited computational resources, deploying MPC remains challenging for strongly nonlinear, high dimensional, and stiff systems \cite{Meadows1996}.
Furthermore, given recent interest in MPC for distributed or mobile systems \cite{Scattolini2009,Christofides2013}, avoiding high computational overhead is even more important for applications where computational resources are limited at the point of control action.
These concerns motivate methods for reducing on-line computational cost in MPC.

To address the demand for excessive computational resources, there has been significant research into ``fast'' MPC over the past two decades.
Broadly speaking, two approaches have been investigated for fast MPC: (1) suboptimal MPC, where a simpler optimization problem is solved that is equivalent (or approximately equivalent) to  the solution of the full problem \cite{Zeilinger2014,Wolf2016FastReview}, and (2) ``explicit'' MPC, where an explicit control law is found that  approximates the implicit MPC controller but does not require on-line optimization \cite{Bemporad2002TheSystems}.
In this work, we focus on the latter, explicit MPC.

For constrained linear systems, an explicit solution can be found by solving a single multiparametric quadratic programming problem off-line \cite{Bemporad2002TheSystems,Tndel2003,Oberdieck2017,Diangelakis2018OnProblems}.
The ``parameters'' are the system states and the solution to this optimization problem is piecewise affine as long as the constraints are linear; with this solution in hand, the controller only needs to first follow a look-up table to determine the relevant polytopic region of state-space in which the system currently lies and then perform an affine computation. 
Unfortunately, for high dimensional systems, the number of polytopes increases exponentially, which can make even the look-up operation too slow for some applications \cite{Chen2018ApproximatingNetworks}.
Multiparametric programming has also been applied to nonlinear MPC, but finding exact solutions to multiparametric nonlinear programming problems is not always feasible, and even when it is feasible, the solution is not necessarily piecewise affine.
Therefore, approximation methods are typically used instead \cite{Dominguez2010RecentProgramming}.

Although the multiparametric programming approach is dominant in the literature, another approach to explicit MPC is interpolation or function approximation.
In this framework, a large number of control policies are computed off-line and the on-line control law is constructed by interpolation (or regression) from states to the optimal policies that were computed off-line.
This method has been most successful using artificial neural networks (ANNs) as the interpolation functions (e.g., \cite{Akesson2006, Karg2018, Chen2018ApproximatingNetworks, Rawlings2019}).
Recently, this method was refined for constrained linear MPC problems using deep neural nets combined with Dykstra's projection to ensure constraint satisfaction \cite{Chen2018ApproximatingNetworks}.

In this work, we present an alternative framework for explicit MPC based on finding an appropriate function between the system state and the control policy that approximates the MPC controller.
We build upon recent advances in nonlinear data mining, especially manifold learning, to show how we can find an intrinsic parametrization of the MPC problem that respects similarities in control policy space.
By working in this intrinsic space, we can design effective interpolating or approximating functions as our explicit controllers.
Specifically, we use state feedback to project the system's position in state space onto a latent manifold, with a parameterization informed by the control policy space and by the state space; and then, using our position in the latent space, we can estimate the entire optimal control policy.
We also demonstrate how we can, in some cases, approximate the ``inverse problem'' of finding the system state \textit{given} the optimal control policy.
Solving the inverse problem could be valuable in contexts outside traditional MPC, where system feedback can be easily observed, but measuring the system state is challenging, as in some biological/living systems.
In order to effectively interpolate between the state space, latent manifold, and control policy, we apply three tools for learning the transformations: polynomial regression (PR), artificial neural networks (ANNs), and Gaussian process (GP) regression.

The remainder of this paper is organized as follows: in Section \ref{sec:theory}, we discuss our approach to MPC using manifold learning and function approximation in more detail; in Section \ref{sec:results}, we present illustrative examples of our approach to several problems including constrained linear systems and a mechanistic model of a nonisothermal continuous stirred tank reactor (CSTR); and in Section \ref{sec:conclusions} we mention  possible applications of this framework and provide some directions for future research.

\section{Theory}\label{sec:theory}
\subsection{Model Predictive Control Background}
First, consider the discrete time nonlinear system written in state space form as a difference equation:

\begin{equation}\label{eqn:gensys} 
\begin{aligned}
x_{k+1} &= f(x_k,u_k) \\
y_k &= h(x_k,u_k) 
\end{aligned}
\end{equation}
where $x\in \mathcal{X}$ is the system state vector, $u\in \mathcal{U}$ is the input vector, $y \in \mathcal{Y}$ is the output vector, and $k \in \mathbb{N}$ is the sample index.
For system \ref{eqn:gensys}, the MPC controller is defined as the controller that minimizes the cost \cite{Mayne2000}:

\begin{equation}\label{eqn:cost}
V(y,k,\mathbf{u})=F(y_{k+N}) + \sum_{i=k}^{k+N-1}\ell(y_i,u_i)
\end{equation}
where $\mathbf{u}=\{u_k,u_{k+1},...,u_{k+N-1}\}$, $\ell(y_i,u_i)$ is the stage cost (such as a normed difference between a reference trajectory, $r_k$, and the predicted trajectory) and $F(y_{k+N})$ is the terminal cost.\footnote{We use upper case, \textbf{Boldface} notation to refer to matrices, lower case, \textbf{boldface} notation to refer to a vector of points in time (e.g., an optimal control policy), and nonboldface notation to refer to vector variables at a single time step,  as well as for scalars.}
$y_i$ for $i>k$ is found by evolving the model (Equation \ref{eqn:gensys}) in time using the policy $\mathbf{u}$ and initial condition $x_k$.
There are often constraints on the range of inputs, range of states, or rate of change of inputs; refer to \cite{Mayne2000} for a full discussion of constrained MPC.
Here, we will assume that we can deploy an optimization algorithm to find the optimal control policy, $\mathbf{u}^*$, that minimizes Equation \ref{eqn:cost}.
Because we are solving the optimization problem off-line, this optimization algorithm need not be highly efficient.
In this paper, we will also assume that the reference trajectory is constant over the entire prediction horizon, i.e., $r_k = r_0$, and could be piecewise constant during process operation; we allow ourselves this strong assumption since the focus is on showcasing the data-driven features of the approach.

The solution to the optimization problem defines a feedback control law.
If many optimization problems are solved off-line, then it is possible, subject to some weak smoothness conditions, to interpolate between the current state, $x_k$, and the optimal control policy, $\mathbf{u}^*$.
The interpolation function is a surrogate for the implicitly defined MPC controller that allows for near instantaneous calculation of the optimal policy.
For first order (i.e., single state), single-input, single-output systems, the problem is trivial; but for more challenging problems, interpolation may be hindered by the ``curse of dimensionality,'' where function training and evaluation time can increase exponentially with dimensionality.
In this work, we discuss recent advances in machine learning and nonlinear data mining that allow us to use tools for function approximation in high-dimensional ({\em but effectively low-dimensional}) spaces, thus enabling explicit computation of the MPC control policy.

\subsection{Manifold Learning of MPC Systems}

Manifold learning algorithms comprise a class of unsupervised machine learning techniques that attempt to find a low-dimensional manifold on which high-dimensional data points are embedded.
The earliest manifold learning algorithm to be developed was principal component analysis (PCA), first presented by Pearson in 1901 \cite{Pearson1901OnSpace}. 
PCA and related techniques, however, are limited in that they can only find a linear embedding of the data. 
Since 2000, a number of related nonlinear manifold learning techniques have been presented \cite{Tenenbaum2000AReduction,Roweis2000NonlinearEmbedding,doi:10.1162/089976603321780317,Coifman2005,Coifman2006,VanDerMaaten2008}.
These methods all share the attribute that they can find nonlinear embeddings; therefore, for example, if 2D data points were all to lie on a smooth curve, PCA would find variation in 2 axes, whereas these manifold learning techniques would (correctly) discover that the data can be represented using only a single variable.

Our key observation is that the optimal control policy, $\mathbf{u}^*$, which has extrinsic dimensionality of $(N\times\dim(u))$ (where $N$ is the number of steps in the control horizon and $\dim(u)$ is the number of manipulated inputs) \textit{must} lie on a manifold with intrinsic dimensionality limited by:

\begin{equation}\label{eqn:dims}
\dim_i(\mathbf{u}^*) \le \dim(x)+\dim(r_0)
\end{equation}
where $\dim(\cdot)$ is the extrinsic dimensionality and $\dim_i(\cdot)$ is the intrinsic dimensionality of an underlying manifold/vector.
This limit results from observing that $\mathbf{u}^*$ is a function of the current system state, $x$, and the current reference trajectory, $r_0$ (the latter, in our case, is just the set point, having dimensionality equal to the number of controlled output variables).

Because of this property, we will refer to the \textit{augmented} state vector, $x^*\in \mathcal{X^*}$, as the concatenation of the state variables and the variables parametrizing the reference trajectory.
The intrinsic dimension may be lower than this maximum if, for example, the system  quickly relaxes to a slow manifold (i.e., it is singularly perturbed), if the  state space realization is nonminimal (i.e., containing redundant information), or if the control policy is a function of the difference between a state variable and a reference variable (as could be the case for linear control of linear systems).

For many systems, however, even though $\mathbf{u}^*$ lies on a manifold of equivalent or lower dimension than $\mathcal{X^*}$, the augmented state space may be very poorly parametrized for predicting $\mathbf{u}^*$.
In other words, a function $c(x^*):\mathcal{X^*}\to\mathcal{U}$ is likely to be quite complicated, requiring many basis functions to represent, and therefore challenging to learn and expensive to evaluate.
Furthermore, it is possible that some (low dimensional) set of nonlinear combinations of $x^*$ are as effective, or nearly as effective, as $x^*$ itself for predicting $\mathbf{u}^*$.
We therefore seek an effective reparametrization of $\mathcal{X^*}$ that prioritizes similarities in policy space, enabling us to find simpler, potentially lower dimensional, representations of the relationship between (augmented) system states and optimal control policies.

In this work, we use the diffusion maps (DMAPS) algorithm for manifold learning \cite{Coifman2005,Coifman2006}, which is reviewed in Appendix \ref{app:dmaps}.
For the MPC problem, we seek to reparametrize the augmented state space in such a way that we can predict the control policy based on knowledge of the system state, as well as the state variables based on knowledge of the control policy.
We will therefore use an ``informed'' metric (c.f., \cite{Holiday2018ManifoldReduction,Lafon2004DiffusionHarmonics}) that considers points ($z$) in an ``input space'' ($\mathcal{Z}$) as well as points ($f(z)$) in a ``function space'' ($\mathcal{F}$) to build the kernel matrix (Equation \ref{eqn:weights}):

\begin{equation}\label{eqn:informed}
d_{i,j} = \frac{\norm{z_i - z_j}^2}{\varepsilon}+\frac{\norm{f(z_i) - f(z_j)}^2}{\xi}
\end{equation}
where $\norm{\cdot}$ is the Euclidean norm, and $\varepsilon$ and $\xi$ are tuning parameters.
We choose the tuning parameters to prioritize distances in the function space, i.e., $\mathrm{median}\left(\norm{z_i - z_j}^2\right) / \varepsilon \ll \mathrm{median}\left(\norm{f(z_i) - f(z_j)}^2\right) / \xi$.
Therefore, it is only at positions where the distance in function space tends to zero that the distance in input space becomes significant.

For the MPC problem, we construct the weight matrix using Equation \ref{eqn:informed} with the augmented state variables used for $z$ and the control policies used for $f(z)$.
We comment that it is not necessary to insert the \emph{entire} control policy ($\mathbf{u}^*$) as $f(z)$ in the distance metric, and instead recommend using $(2\dim(x^*)+1)$ steps of the control policy (as a heuristic that was inspired by the Takens embedding theorem \cite{Takens1981}).
By appropriate choice of $\varepsilon$ and $\xi$, the informed metric finds a parametrization of the augmented state space that is organized primarily by similarities in policy space---when $\norm{\mathbf{u}^*_i - \mathbf{u}^*_j}$ is large, it will spread these points far apart in the intrinsic space.
Thus, the intrinsic parametrization will organize the control policy space and the state variable space using as few variables as possible---subject to the limit in Equation \ref{eqn:dims}---and hopefully make prediction of the high-dimensional control policies ``simpler''.
The leading nonredundant eigenvectors can be used to parametrize the manifold of the outputs (the control policies) and therefore provide the coordinates important for predicting them \cite{Coifman2005, Coifman2006}.
Therefore (nonredundant \cite{Dsilva2015}) eigenvectors---sometimes after a large gap in the spectrum---correspond to coordinates in the augmented state space that are unimportant for predicting the control policy and can be eliminated to provide a reduced order approximation of the control law.

\begin{figure}
\begin{center}
\includegraphics[height=3cm]{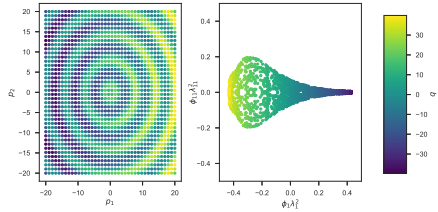}
\caption{(a) Points sampled from $\mathbb{R}^2$ colored by a ``complicated'' function. (b) Intrinsic reparametrization of $\mathbb{R}^2$ learned by DMAPS using the informed metric in Equation \ref{eqn:informed} (see text for discussion).}
\label{fig:informed_demo}
\end{center}
\end{figure}

We will call the manifold discovered by this procedure ``the control manifold'', $\mathcal{M}_c$; the non-redundant DMAPS eigenvectors , $\{\phi\lambda^k\}$, provide its parametrization (where the eigenvectors are optionally scaled by their corresponding eigenvalues, or powers of their eigenvalues, see Appendix \ref{app:dmaps}). 
We note that $\mathcal{M}_c$ has dimensionality less than or equal to that of $\mathcal{X}^*$, and typically much less than the extrinsic of the control policies.

Because the parametrization of $\mathcal{M}_c$ respects similarities in policy space, we expect that predicting $\mathbf{u}^*$ given $\phi$ should be a ``simple'' task.
If $\norm{\mathbf{u}^*_i - \mathbf{u}^*_j}$ is small, then $\norm{\phi_i - \phi_j}$ should be small as well (unless $\norm{x^*_i-x^*_j}$ is very large).
To demonstrate, consider the function $q=f(p):\mathbb{R}^2\to \mathbb{R}$.
Here, $q=f(p) = 10\sin{\sqrt{p_1^2+p_2^2}}+p_2$.
Data ($p$) sampled from a regular grid in $\mathbb{R}^2$ and colored by function value ($q$) are shown in Figure \ref{fig:informed_demo}a.
If we did not \textit{a priori} know the function, we may resort to complicated and/or computationally expensive methods to approximate it in ($p_1, p_2$) space.
However, in the reorganized space discovered via DMAPS and shown in Figure \ref{fig:informed_demo}b, the values $q$ are seen to be a  ``simple'' function of the single DMAPS coordinate $\phi_1$.

\subsection{Function Approximation}
The above discussion implies that it is possible to link $x^*$ to $\mathbf{u}^*$ through $\mathcal{M}_c$, but does not discuss how to find our position on $\mathcal{M}_c$ in $\phi$ coordinates, nor, how, given $\phi$, to predict the policy, $\mathbf{u}^*$.
Additionally, we would like to solve the inverse problem, which requires predicting $\phi$ from $\mathbf{u}^*$, and then predicting $x^*$ from $\phi$, with the important caveat that the problem is not in general invertible, primarily when the control policy pushes against constraints or when the system dynamics are overspecified.
The DMAPS metric from Equation \ref{eqn:informed} was designed to make the task of predicting $\mathbf{u}^*$ from $\phi$ (or vice versa) simple in comparison with predicting $\mathbf{u}^*$ from $x^*$, since the former function is defined over the lower-dimensional, intrinsic space.
In any case, we use three methods for function approximation, each of which has advantages and disadvantages: polynomial regression (PR), artificial neural networks (ANNs), and Gaussian process (GP) regression.
Other tools for function approximation, such as geometric harmonics \cite{Coifman2006GeometricFunctions} or Laplacian pyramids \cite{Rabin2012HeterogeneousPyramids} may be used as well.
In principle, any of these methods may be applicable to learn any of the functions of interest ($x^* \to \phi$, $\phi \to \mathbf{u}^*$, and their inverses), though in practice PR is unlikely to effectively capture the mapping between $x^*$ and $\phi$, which is typically quite complicated.

First, we examine polynomial regression, a classical linear technique.
In some cases, the mapping between $\phi$ and $\mathbf{u}^*$ is sufficiently simple that this classical approach works quite well (see Figure \ref{fig:informed_demo}, for example).
We use the ordinary least squares estimator, $\hat{\Theta}=(X^TX)^{-1}X^TY$ to find the coefficients of the polynomial regression estimator, where $X$ is a feature matrix and $Y$ is an output matrix.
For this problem, rows of $X=\begin{bmatrix} \phi & \phi^2 & ...& \phi^n\end{bmatrix}$ and rows of $Y=\mathbf{u}^*$.
Then, the output at new values of $\phi$ can be predicted using $\hat{Y}=\hat{\Theta} X$ \cite{hastie01statisticallearning}.

Artificial neural networks have become popular in many applications due to their versatility and ability to represent arbitrary continuous functions, and are widely considered the workhorse method for ``deep learning'' \cite{Nielsen2015,Lecun2015}.
Artificial neural networks have also been widely used in control for decades, mostly for empirical approximation of state equations \cite{HUNT19921083}, but also for explicit MPC problems \cite{Chen2018ApproximatingNetworks, Akesson2006, Karg2018, Rawlings2019}.
The main disadvantage of ANNs is that they require significant computational time for training, and require manual design of network topology and choice of activation function.
More complex networks, with more nodes and hidden layers, require more training time and are susceptible to overfitting.

Finally, GP regression, as a non-parametric Bayesian modeling technique, provides the conditional distribution of an output as a function of its observed inputs.
In order to deal with higher dimensional input spaces, we introduce the automatic relevance determination (ARD) weight in the covariance kernel, which employs an individual lengthscale hyperparameter for each input dimension \cite{Rasmussen:2005:GPM:1162254}.
In this paper, we employ a Mat\'{e}rn kernel for the covariance:

\begin{equation}\label{eqn:matern}
\mathbf{\kappa}(x_i,x_j) = \left(1+\sqrt{3}\theta\mathrm{d}(x_i,x_j)\right)\exp \left( -\sqrt{3}\theta\mathrm{d}(x_i,x_j)\right).
\end{equation}

GP regression typically requires some training time, in that hyperparameters should be optimized; but there are usually far fewer adjustable parameters than in ANNs.

\section{Results and Discussion}\label{sec:results}
To illustrate our approach, we present (1) for validation purposes, an instructive demonstration of the paramatrization learned by diffusion maps for constrained linear systems and (2) a tutorial example of our DMAPS-enabled explicit MPC formulation for controlling a jacketed nonisothermal continuous stirred tank reactor (CSTR).

\subsection{Constrained MPC for Linear Systems}
Nearly all practical MPC problems have constraints on the state variables, input variables, or rates of change of input variables.
Therefore, even in the case of a linear system model, the control law is nonlinear.
As long as the constraints are linear, the control law will be  piecewise affine with a set of so-called ``critical regions'' in the state space that correspond to which constraints are active \cite{Bemporad2002TheSystems}.
We will use our manifold learning techniques to learn effective parametrizations for these kinds of nonlinear control laws.

First, to illustrate how the DMAPS algorithm will automatically detect an overspecified problem and identify a lower-dimensional parametrization, we introduce the regulatory control problem for the following singularly perturbed system:
\begin{equation} \label{eqn:sing_perturb}
x_{k+1} =\begin{bmatrix} 0.4079 & 0.4031 \\
                        0.4157 & 0.4109\end{bmatrix}x_k + 
        \begin{bmatrix} 0.7071 \\ 0.7071 \end{bmatrix}u_k  
\end{equation}
where we assume exact state information is available ($y_k=x_k$). Note that the eigenvalues of the state transition matrix are $\{0.819, 4.54E-5\}$, indicating significant time scale separation that results in an essentially one-dimensional system; furthermore, the Hankel singular values (which characterize input to output energy transfer for each state when the system is in balanced coordinates \cite{Moore1981}) are $\{3.03, 7.81E-3\}$, again suggesting the system is effectively one-dimensional.
For this problem, the MPC control policy can be found via the solution to the following constrained quadratic programming problem:

\begin{equation}\label{eqn:MPCsing_perturb}
\begin{aligned}
& \underset{\mathbf{u}} {\text{min}} 
& V(x,k,\mathbf{u}) &= x_{k+5}^T\mathbf{Q_t}x_{k+5} + \sum_{i=k}^{k+4}\left({x_i^T\mathbf{Q}x_i + u_i^T\mathbf{R}u_i}\right)\\
&\text{s.t.} & u_i &\le 0.5,~i=k, k+1, ..., k+4 \\ 
& & -u_i &\le 0.5,~i=k, k+1, ..., k+4 \\
\end{aligned}
\end{equation}
where $\mathbf{Q_t}$ is the terminal cost (found by solving the discrete time Lyapunov equation $\mathbf{Q_t}=\mathbf{A}^T\mathbf{Q_t}\mathbf{A}+\mathbf{Q}$; $\mathbf{A}$ is the state transition matrix for System \ref{eqn:sing_perturb}), $\mathbf{Q}=\mathbf{I}$ is the state cost, and $\mathbf{R}=0.01$ is the control cost.

We sample a 20 x 20 regular grid from the state space, and calculate the control policy at each of the grid points.
As recommended in Section \ref{sec:theory}, we insert the state variables for $z$ and the first 5 steps of the control policy as $f(z)$ into Equation $\ref{eqn:informed}$, and then construct the DMAPS embedding.
Additionally, we use PCA (where each row in the data matrix is the concatenation of the state vector and the control policy), to compare the nonlinear manifold learning technique with a classical linear method.

\begin{figure}
\begin{center}
\includegraphics[height=6cm]{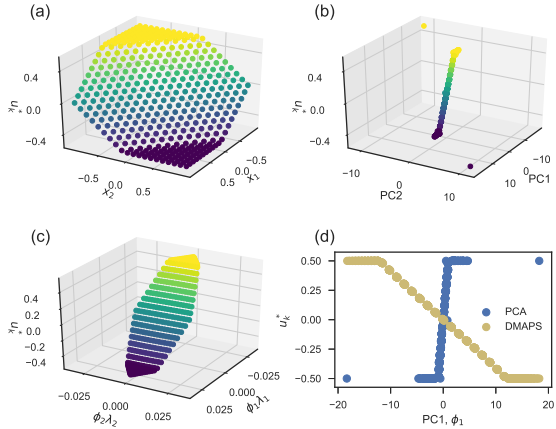}
\caption{MPC control law from the the solution to the optimization problem in Equation \ref{eqn:MPCsing_perturb}, which is a singularly perturbed linear problem with linear constraints, plotted in several coordinate systems: (a) The grid of points sampled from the state space, (b) the first two principal components, and (c) the DMAPS embedding coordinates. (d) The control action plotted against the first PC or DMAPS eigenvector, demonstrating how a low dimensional control law is discovered (rescaled for easier visualization).}
\label{fig:lin_sing}
\end{center}
\end{figure}

The MPC control law (i.e., the first step of the control policy) is shown as a function of the states in Figure \ref{fig:lin_sing}a and, as anticipated, is effectively one dimensional (along the diagonal).
Notice that there are apparently three critical regions for predicting the control law: the minimum and maximum input constraints, and a plane that connects them.
For this simple problem, DMAPS and PCA can both find effective 1D embedding spaces (the first coordinate in Figures \ref{fig:lin_sing}b-c), and an exact 1D control law (Figure \ref{fig:lin_sing}d).

Next, we investigate the system first examined by Bemporad et al. \cite{Bemporad2002TheSystems} and demonstrate how we can learn the relevant control laws for the full-dimensional system, as well as find an effective reduced order approximation to the control law. Bemporad et al, introduced the system:

\begin{equation} \label{eqn:bempsys}
x_{k+1} =\begin{bmatrix} 0.7326 & -0.0861 \\
                        0.1722 & 0.9909\end{bmatrix}x_k + 
        \begin{bmatrix} 0.0609 \\ 0.0064\end{bmatrix}u_k  
\end{equation}
where again we assume that we have perfect state knowledge.
The system in Equation \ref{eqn:bempsys} has the Hankel singular values  $\begin{bmatrix}0.4445 & 0.1522\end{bmatrix}$, suggesting that approximate model order reduction may be possible.

First, consider the regulatory control problem of driving System \ref{eqn:bempsys} to the origin, subject to the constraint that $-2<u_i<2$.
The following constrained quadatric programming problem can be solved to provide the MPC control policy:

\begin{equation}\label{eqn:MPClin1}
\begin{aligned}
& \underset{\mathbf{u}} {\text{min}} 
& V(x,k,\mathbf{u}) &= x_{k+10}^T\mathbf{Q_t}x_{k+10} + \sum_{i=k}^{k+9}\left({x_i^T\mathbf{Q}x_i + u_i^T\mathbf{R}u_i}\right)\\
&\text{s.t.} & u_i &\le 2.0,~i=k, k+1, ..., k+9 \\ 
& & -u_i &\le 2.0,~i=k, k+1, ..., k+9 \\
\end{aligned}
\end{equation}
We find the terminal cost using the same procedure as in the previous problem and, as before, $\mathbf{Q}=\mathbf{I}$ and $\mathbf{R}=0.01$.
Note that this is the same control problem studied by Bemporad \cite{Bemporad2002TheSystems}, except that we expand the time horizon from 5 to 10 time steps.

\begin{figure}
\begin{center}
\includegraphics[height=6cm]{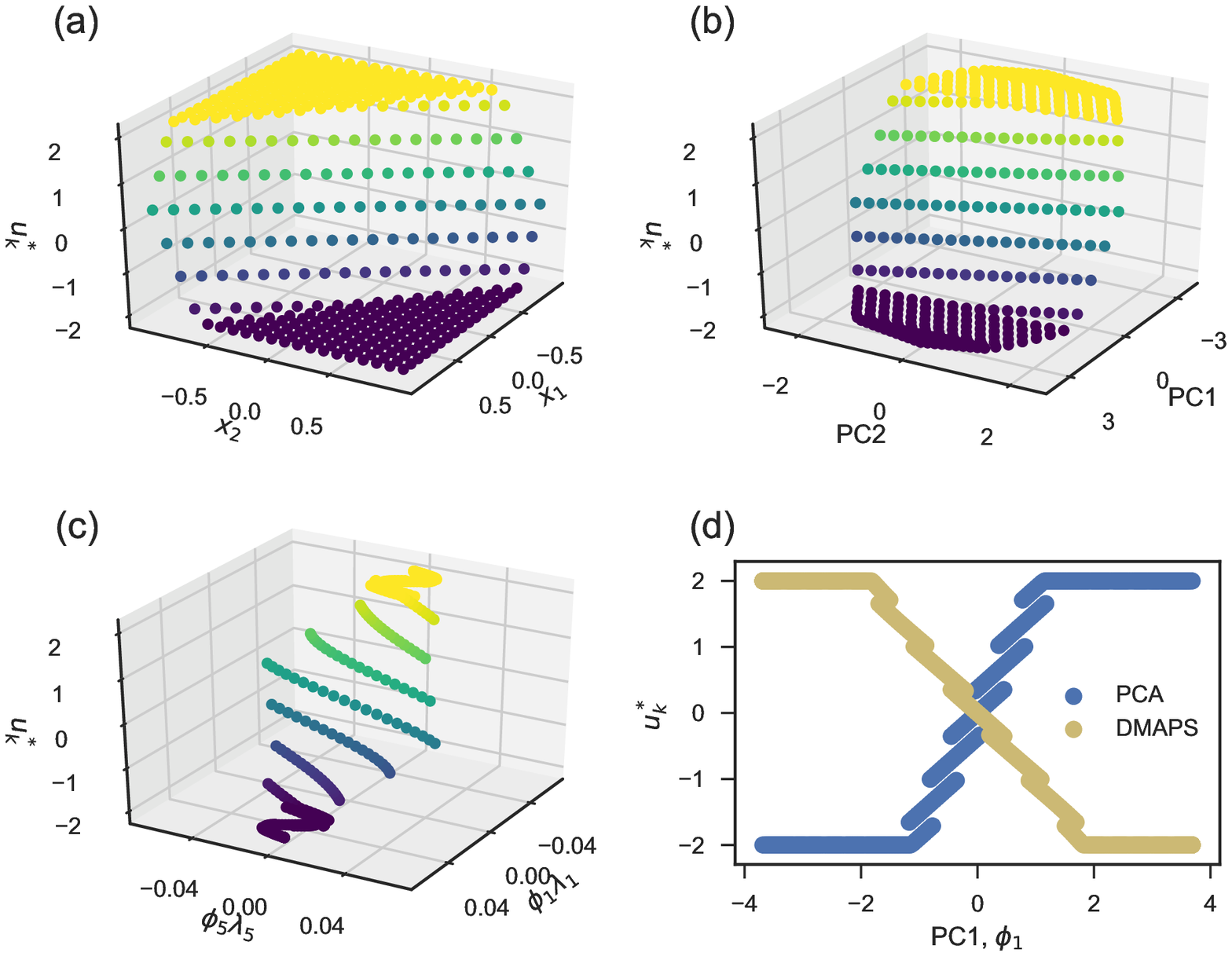}
\caption{MPC control law from the the solution to the optimization problem in Equation \ref{eqn:MPClin1}, which is a linear problem with linear constraints, plotted in several coordinate systems: (a) The grid of points sampled from the state space, (b) the first two principal components, and (c) the DMAPS embedding (d) The control action plotted against the first PC or DMAPS eigenvector, demonstrating the existance of an \emph{approximate} one-dimensional control law (rescaled for easier visualization).}
\label{fig:lin_con1}
\end{center}
\end{figure}

We follow the same procedure as for the singularly perturbed problem to generate the control law in the original space, a PCA embeddeding space, and a DMAPS embedding space, as shown in Figure $\ref{fig:lin_con1}$.
Here, we note that while once again both PCA and DMAPS find effective 2D parametrizations of the control law (Figure \ref{fig:lin_con1}b-c), and approximate reduced order control laws (Figure \ref{fig:lin_con1}d), the reduced order DMAPS control law is noticeably improved in comparison to PCA.
The improvement is both because DMAPS can uncover nonlinear relationships among the variables and because we can design the kernel in a way that favors distances in control policy space. 

The previous example demonstrated how DMAPS can reorganize, and reduce, that state dimension for predicting the control action.
Now, we return to the same linear system (Equation \ref{eqn:bempsys}), but following Diangelikis et al, \cite{Diangelakis2018OnProblems} (who found the explicit control law using multiparametric quadratic programming), we introduce quadratic constraints.
The optimization problem defining the control policy is now:
\begin{equation}\label{eqn:mpcquad}
\begin{aligned}
& \underset{\mathbf{u}} {\text{min}} 
& V(x,k,\mathbf{u}) &= x_{k+6}^T\mathbf{Q_t}x_{k+6} + \sum_{i=k}^{k+5}\left({x_i^T\mathbf{Q}x_i + u_i^T\mathbf{R}u_i}\right)\\
&\text{s.t.} & u_i &\le 2.0,~i=k, k+1, ..., k+5 \\ 
& & -u_i &\le 2.0,~i=k, k+1, ..., k+5 \\
& & u_k^2 &\le x_k^T x_k
\end{aligned}
\end{equation}
where we have added the constraint $u_k^2\le x_k^T x_k$. $\mathbf{Q_t}, 
\mathbf{Q}$ and $\mathbf{R}$ are found using the same approach as before; we reduced the control horizon to 6 steps (for computational efficiency).

\begin{figure}
\begin{center}
\includegraphics[height=6cm]{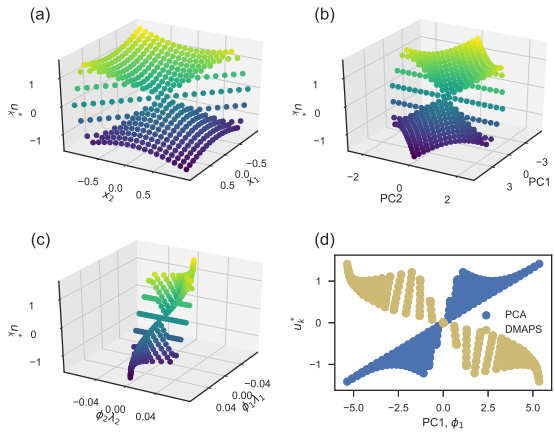}
\caption{MPC control law from the the solution to the optimization problem in Equation \ref{eqn:mpcquad}, which is a linear problem with quadratic constraints, plotted in several coordinate systems: (a) The grid of points sampled from the state space, (b) the first two principal components, and (c) the DMAPS embedding. (d) The control action plotted against the first PC or DMAPS eigenvector, indicating that the control action \emph{cannot} be written as of function of the first PC or DMAPS eigenvector, and therefore a reduced order form would be a poor approximation (rescaled for visualization)}
\label{fig:lin_quadcon}
\end{center}
\end{figure}

Again, we generate the control policy on a 20 x 20 regular grid, and find DMAPS and PCA embeddings of the results using the same procedure as before.
Now, however, we observe in Figure \ref{fig:lin_quadcon} that nearly exact order reduction is no longer possible.
There are still three apparent critical regions seen from plotting $u_k^*$ as a function of $x$ (Figure \ref{fig:lin_quadcon}a), but now two of them are paraboloids rather than planes.
In any case, the first two principal components and the first two DMAPS eigenvectors provide an appropriate parametrization (Figure \ref{fig:lin_quadcon}b-c).
Because the control law is intrinsically 2D (even though the dynamics are approximately 1D), DMAPS cannot find an effective reduced order representation; and, of course, neither can PCA (Figure \ref{fig:lin_quadcon}d).

The previous examples illustrate how DMAPS finds an intrinsic parametrization of the system state that is organized by the control policy.
We demonstrated how, when possible, DMAPS will suggest a lower dimensional manifold from which we can write the control law.
In the next section, we will use these same tools as presented here for a more complicated, nonlinear system, and further demonstrate how to use function approximation techniques to write the control law in an explicit form.

\subsection{Nonlinear MPC for a Nonisothermal CSTR}
\begin{figure}
\begin{center}
\includegraphics[height=3cm]{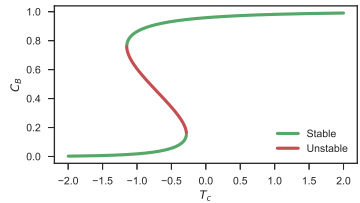}
\caption{Bifurcation diagram showing steady-state concentration of species $B$ (the controlled variable) as a function of the cooling water temperature (the manipulated variable).}
\label{fig:bif}
\end{center}
\end{figure}

We demonstrate the use of our manifold learning approach to approximate the MPC control law for a nonlinear, unstable chemical reactor system where both reference tracking and disturbance rejection is considered.
The CSTR contains a single reaction, $A\to B$, which is exothermic with reaction rate temperature dependence having the Arrhenius form.
We will control the concentration of species $B$ using the temperature of the cooling water as the manipulated variable.
The cooling water temperature is constrained both with minimum and maximum absolute values and with a maximum stepwise rate of change.
All of the data presented are expressed in dimensionless variables.
For detailed descriptions of the model equations and the MPC parameters, refer to Appendix \ref{app:cstr} and Appendix \ref{app:mpc}, respectively.

\begin{figure}
\begin{center}
\includegraphics[height=6cm]{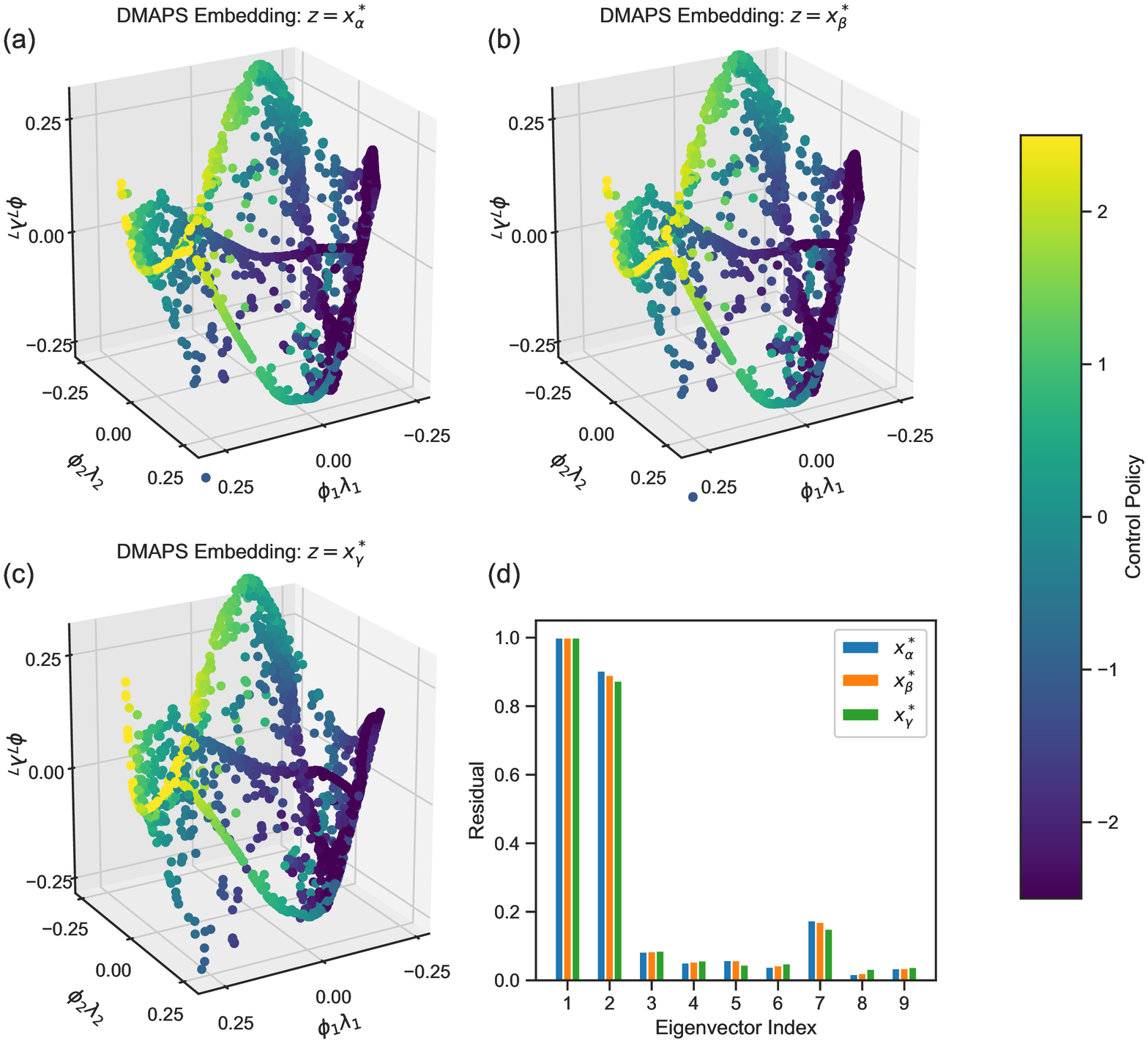}
\caption{DMAPS embedding of the MPC systems showing the training data using three alternative state parametrizations (a) $x_\alpha^*$, (b) $x_\beta^*$, and (c) $x_\gamma^*$ and the first 10 steps of the control policy ($\mathbf{u}^*$) in the distance metric (see text). (d) Residuals from local linear regression for each parametrization in (a)-(c) suggesting that the first, second, and seventh eigenvectors are the best parametrization of the underlying manifold (see Equation \ref{eqn:LLR} and Reference \cite{Dsilva2015}). Arguably, only the first two eigenvectors are necessary, though we observe improved prediction accuracy when including the next nonredundant eigenvector.}
\label{fig:alt_ss}
\end{center}
\end{figure}

\begin{figure*}
\begin{center}
\includegraphics[height=6cm]{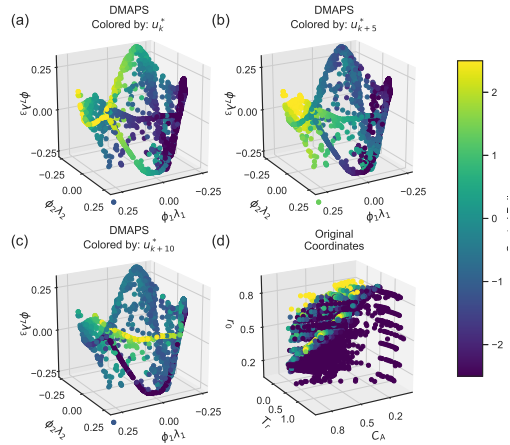}
\caption{Training data using $x_\alpha^*$: (a-c) The first three DMAPS intrinsic coordinates learned using the informed metric and colored by the first, fifth, and tenth steps of the control policy (d) The original augmented states, $x_\alpha^*$, colored by $u_k^*$ (the first step of the control policy).}
\label{fig:dmaps}
\end{center}
\end{figure*}

\begin{figure*}
\begin{center}
\includegraphics[height=6cm]{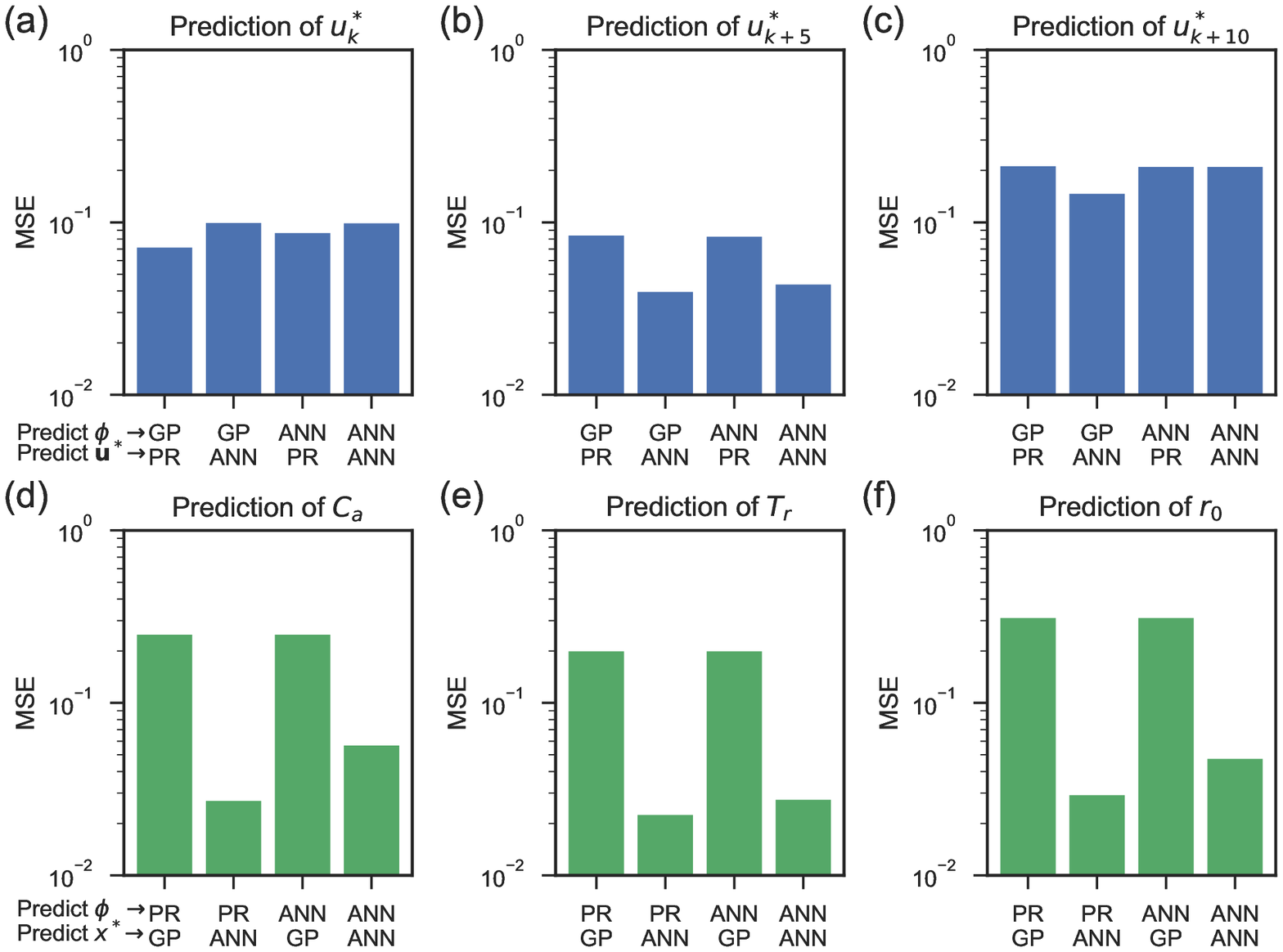}
\caption{Model performance using test data: (a-c) Mean squared prediction error for different explicit approximations of the first, fifth, and tenth steps of the control policy using various function approximation techniques. (d-f) Mean squared prediction error for the inverse problem (finding the augmented states from the control policy), using various function approximation techniques.}
\label{fig:validation}
\end{center}
\end{figure*}

The nonisothermal CSTR is challenging to control partly because a large range of state space is unstable at steady state.
Figure \ref{fig:bif} shows the bifurcation diagram of the open loop system with cooling water temperature as the bifurcation parameter.
The system contains two saddle-node bifurcations and requires feedback stabilization to operate at conversions between approximately 0.2 and 0.8.
If using MPC, this means that short sampling intervals are required to ensure that feedback is frequent enough to stabilize the system.
Because the sampling interval is short, that likewise means that the control horizon must be long, so that significant dynamics are included in the cost function.
These challenges make this system an appropriate candidate to test our formulation of explicit MPC control.

One of the key advantages of the manifold learning approach is that it will automatically detect equivalent or redundant descriptions of the state space in a purely data-driven manner.
For example, we know from theory that our process model is intrinsically second order and can be described using the physical variables of $C_A$ and $T$ (see Appendix \ref{app:cstr}).
For more complicated processes, however, we may not recognize the intrinsic state space and instead model our process using whichever measured variables are at hand to represent it.
For a purely ANN-driven approach, using such alternative (and possibly redundant) parametrizations may require redesigning the network topology for the new problem, and certainly requires retraining the network.
The manifold learning approach, however, will always detect an intrinsic 2D space, regardless of how the state space was parametrized.

To illustrate, we consider three alternative parametrizations of state space: (1) $x_\alpha=\begin{bmatrix} C_a & T_r\end{bmatrix}$, the original variables used for modeling the system, (2) $x_{\beta}=\begin{bmatrix}\left(kC_a\right) & \left(\frac{q}{V}(T_0-T_r)- \frac{\Delta H}{\rho C_p}C_A\right)\end{bmatrix}$, the reaction rate and heating rate (see Appendix \ref{app:cstr}) and (3) $x_\gamma=\begin{bmatrix}x_\alpha & x_\beta\end{bmatrix}$, the concatenation of the other two parametrizations.
By concatenating the reference variable to each of these, we define three augmented state vectors, $x_\alpha^*$, $x_\beta^*$ and $x_\gamma^*$, any one of which is a complete description of the system state and sufficient (in principle) to predict the control policy.

To discover $\mathcal{M}_c$, the control manifold, we need to sample the augmented state space and compute the control policies off-line.
We sampled 200 points randomly in augmented state space, and evolved the system in time for 20 time steps using a nonlinear model predictive controller with time horizon of 20 time steps.
Taken together, we have 4000 samples, and we collected all of the control policies in matrix $\mathbf{U}^*\in \mathbb{R}^{4000 \times 20}$ and augmented state variables in matrices $\mathbf{X}_\alpha^* \in\mathbb{R}^{4000 \times 3}$, $\mathbf{X}_\beta^* \in\mathbb{R}^{4000 \times 3}$, and $\mathbf{X}_\gamma^* \in\mathbb{R}^{4000 \times 5}$.
To test our tools for function approximation, we randomly partitioned the data into 3000 points for training and 1000 points for testing.

We apply DMAPS three times using the informed metric from Equation \ref{eqn:informed} to the downsampled policy matrix $\mathbf{U}^*$ with each of the augmented state matrices $\mathbf{X}_\alpha^*$, $\mathbf{X}_\beta^*$, and $\mathbf{X}_\gamma^*$. 
The results for each case are shown in Figure \ref{fig:alt_ss}, which shows augmented state variables projected into the intrinsic space and colored by $u_k^*$.
These results demonstrate that DMAPS can find an effective, possibly lower dimensional, intrinsic manifold regardless of which state variables are made explicit, from which the control policy can be predicted.
For each of the different parametrizations, we find that the first three eigenvalues are the best intrinsic parametrization (arguably, based on the residuals in Figure \ref{fig:alt_ss}d, only the first 2 DMAPS eigenvectors are needed---however, we observe in practice that including the third improves prediction accuracy).

Using $x_\alpha^*$ (for example) as the augmented state variables, Figure \ref{fig:dmaps} shows that given the position in DMAPS coordinates, we can predict the control policy---not just the control law, $u_k^*$---using our intrinsic parametrization.
The policy space can be effectively parametrized using only 3 intrinsic dimensions, $\phi=\begin{bmatrix}\phi_1& \phi_2 & \phi_3 \end{bmatrix}$.
Facilitated prediction using this intrinsic space parametrization is possible because the mapping from $\phi$ to $\mathbf{u}^*$ is much simpler than from $x_\alpha^*$ to $\mathbf{u}^*$ (as is visually evident by comparing Figure \ref{fig:dmaps}a and \ref{fig:dmaps}d).
Therefore, simpler functions can be used to determine the control policy when using the latent space rather than the original space.

Now, using $x_\alpha^*$ as inputs, we design explicit feedback control laws, i.e., functions  $c(x_\alpha^*):\mathcal{X}_\alpha^*\to\mathcal{U}$.
This task is divided into two stages: first, estimate the intrinsic variables $\phi$ given the augmented state variables $x_\alpha^*$, and second, predict $u_k^*$ given $\phi$.
We emphasize that because we transformed to the intrinsic variables, we can also easily predict several steps of the control policy, rather than just the first step.

For estimating $\phi$ from $x_\alpha^*$, we use two methods presented in Section \ref{sec:theory}: ANNs and GPs.
We write the neural network model as:

\begin{equation} \label{eqn:nn1}
\hat{\phi}_{ANN} = p(x_\alpha^*;\mathbf{W})
\end{equation}
where $p$ is an artificial neural network with eight hidden layers of 20 nodes each, three input nodes corresponding to each element in $x_\alpha^*$, three output nodes corresponding to each element in $\phi$, and $\mathbf{W}$ is the weight matrix for the neural network.
Rectified linear activation functions (which have become a \textit{de facto} standard in deep learning \cite{Lecun2015}) were used for each of the nodes in the hidden layers:

\begin{equation}\label{eqn:relu}
f(x) = \max(0,x)
\end{equation}
To facilitate the use of the activation function in Equation \ref{eqn:relu} (which will never predict a negative output), all of the data variables were linearly rescaled from 0 to 1, and scaled back to their original ranges for presentation.
Equation \ref{eqn:nn1} was built using pyTorch and trained using the Adam optimizer (a similar algorithm to stochastic gradient descent) with learning rate of $1\times 10^{-3}$ \cite{Kingma2014,Paszke2017}.
Alternatively, we predict $\phi$ from $x_\alpha^*$ using GP regression.
Using the Mat\'ern covariance function in Equation \ref{eqn:matern}, we optimize hyperparameters for prediction over our training data set by minimizing negative log marginal likelihood \cite{Rasmussen:2005:GPM:1162254} to obtain $\hat{\phi}_{GP}$.

With intrinsic variables $\hat{\phi}$ in hand using one of the three methods above, we now predict $\mathbf{u}^*$.
We use cubic PR (higher order polynomial regression was tested, with no improvement) and an ANN model.
As before, the ANN model uses the rectified linear activation function from Equation \ref{eqn:relu} (with inputs and outputs appropriately rescaled), eight hidden layers of 20 nodes each, three input layers for $\phi$, and 20 output neurons for $\mathbf{u}^*$.

Using two of the techniques for predicting $\phi$ and two of the techniques for predicting $\mathbf{u}^*$, we have designed four explicit MPC controllers.
Their prediction accuracy is shown for the first thre steps of the control policy in Figure \ref{fig:validation}a-c.
All of the methods provide comparable prediction accuracy---though we comment that we did not spend effort optimizing neural network performance and may see some improvement with different architectures, activation funcions, or optimization procedures.
Additionally, we did not prioritize learning any particular step of the control policy---if we had weighted the loss function to give more emphasis on (for example) $u_k^*$ than the subsequent steps, we likely would have lower error for $u_k^*$ (which defines the control law). 

We also solve (where feasible) the inverse problem: deducing the augmented state space from the control policy, by following the same procedure as above in reverse.
We find PR and ANN models to predict $\phi$ from $\mathbf{u}^*$, and then find ANN and GP models to predict $x_\alpha^*$ from $\phi$.
For this second task, we note that when the control policy pushes against a constraint (defined here as $u_k^* = -2.5$ or $u_k^* = 2.5$) predicting the augmented state accurately becomes infeasible because the transformation is not invertible.
Figure \ref{fig:validation}d-f shows the mean squared error of our state predictions based on observations of control policies for test points that are not against constraints.

\begin{figure}
\begin{center}
\includegraphics[width=.5\textwidth]{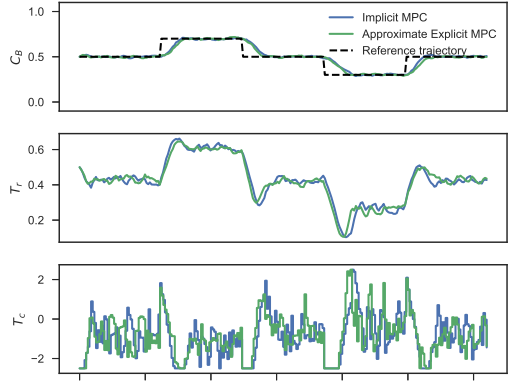}
\caption{Comparison of implicit MPC controller with DMAPS-enabled approximate explicit nonlinear MPC controller using GP regression to predict $\phi$ from $x_\alpha^*$ and an ANN model to predict $u_k^*$ from $\phi$.}
\label{fig:mpc_results}
\end{center}
\end{figure}

Finally, we test the control laws we developed on-line using the full simulated CSTR model.
Figure \ref{fig:mpc_results} shows the performance of our DMAPS enabled, explicit nonlinear model predictive controller, in comparison to the implicit (and therefore exact) MPC controller.
We use state observations, $x_{\alpha,k}^*$, to predict $\hat{\phi}_k$ using GP regression and then an ANN models to predict $\hat{u}_k^*$ at each time step.
The reference trajectory is open-loop unstable (see Figure \ref{fig:bif}), and white Gaussian noise is added to each state as a disturbance.
Still, the controller can maintain the process near the reference value, as well as effectively track set point changes, while closely matching the implicit controller.

\section{Conclusions and Future Directions}\label{sec:conclusions}
In this paper, we have developed and demonstrated a data-driven approach to designing explicit model predictive controllers.
We uncover intrinsic, low-dimensional structure in the high-dimensional control policies and state vectors that provides a link between state space and policy space.
Our manifold learning method is agnostic to the particular parametrization of state space, and equally valid regardless of which state variables are available (as long as we have sufficiently
many).
Furthermore, because the similarity measure used by DMAPS is designed to favor similarities in policy space, predicting the entire control policy becomes about as easy as predicting just its first step.
Like other approaches to explicit MPC, by developing functions between the state space and control action, we avoid the need for on-line optimization.
Although our tutorial demonstrations were for single-input, single-output control problems, our framework naturally generalizes to multiple-input, multiple output systems---in fact, we showed that we can easily generate multiple outputs in that we can predict the full time series of control actions.

We showed how the DMAPS algorithm will identify redundant information in the control policy space, and thus can identify a reduced order control law even when the augmented state space is densely sampled.
The reduced order model can be due to unbalanced coordinates, singular perturbation, or other redundant specifications of the systems.
As illustrated by the quadratically constrained problem, we notice that by directly learning the manifold on which the control policy lies, rather than the manifold on which the dynamics lie, we automatically detect cases where even though the dynamics are low order, the control policy (and therefore, the relevant MPC problem) is not low dimensional due to constraints.

Ensuring constraint satisfaction is a challenge for any approach that relies on function approximation.
For cases where it is feasible, we could use Dykstra's projection, to project the control law into the feasible region, as recommended by Chen et al \cite{Chen2018ApproximatingNetworks}.

Here, we assumed that we have access to the full system state, $x$, either via direct measurement or from a state estimator, at each sampling point.
Inspired by delay embedding theorems, such as that of Takens \cite{Takens1981}, we expect that we do not need full state feedback to apply our methodology.
Much like a state estimator synthesizes information from histories of measurements to ``observe'' the system state, we could directly ``observe'' the control policy using this information.
In future work, we will investigate how to build a purely data-driven ``policy observer,'' where the MPC policy is predicted from only measured data.
We anticipate that, as long as the usual state observability conditions are satisfied, the control policy will be equally observable.

Alternatively, there may be situations where feedback is infrequent, but actuation is comparatively fast.
Here, we could use our framework to design an explicit controller that takes control action in between state observations.
Our procedure is especially well-suited for such systems because we can easily estimate the full, $N$-step control policy, and not only the first step like other explicit controllers that use function approximation.
Finally, even if not accurate enough to implement, our data-inferred control policy can conceivably be used for ``smarter'' initialization of the optimization algorithm in implicit MPC. 

We also demonstrated how, given observations of the control policy, we can predict the system state.
From the perspective of process control, where the objective is developing a feedback control law, the inverse problem may not appear to be of interest.
Yet, we can imagine a system in which the ``feedback'' is easy to observe, while the system state is not.
For example, we may consider an ``expert human'' controller, or a neural network controller.
Using this approach, we could design a state observer that uses control observations to determine that state, instead of output measurements like in conventional state observers.

\section*{Acknowledgements}
The assistance of Dr. Mahdi Kooshkbaghi with ANN training is gratefully acknowledged.

\bibliographystyle{plain}
\bibliography{main}

\appendix

\section{Diffusion Maps} \label{app:dmaps}
DMAPS is an algorithm for manifold learning \cite{Coifman2005,Coifman2006}.
We assume that the data points $X \in \mathbb{R}^N$ all lie on a smooth low dimensional manifold, $\mathcal{M}$, such that $\dim(\mathcal{M})\ll N$.
The algorithm works by finding data-driven approximations for eigenfunctions of the Laplace-Beltrami (i.e., diffusion) operator on  $\mathcal{M}$.
One can show that these eigenfunctions provide an effective minimal parameterization of $\mathcal{M}$ \cite{Jones2008}.

To approximate the eigenfunctions, we find eigenvectors of the appropriately normalized graph Laplacian, which is constructed as follows. Given $m$ observations of the data, define a kernel matrix, $\mathbf{W}\in \mathbb{R}^{m\times m}$:

\begin{equation}\label{eqn:weights}
\mathbf{W}_{ij} = \exp{\left(-d_{i,j}\right)},i,j=1,...,m
\end{equation}
where $d_{i,j}$ is a metric representing distance between data points $i$ and $j$ that is often the Euclidean distance, though in this work we use an input-output informed metric defined in Equation \ref{eqn:informed}.
To account for nonuniform sampling, the kernel matrix is normalized with $P_{ii}=\sum_{k=1}^m W_{ik}$ using:

\begin{equation}
\tilde{\mathbf{W}}=\mathbf{P}^{-\alpha}\mathbf{W}\mathbf{P}^{-\alpha}
\end{equation}
where for isotropic DMAPS, $\alpha=1$.
Next, we define a diagonal matrix from the row sums of the kernel matrix:

\begin{equation}\label{eqn:diag}
\mathbf{D}_{ii}=\sum_j\tilde{\mathbf{W}}_{ij}.
\end{equation}
Finally, we construct the Markov transition matrix:

\begin{equation}\label{eqn:markov}
\mathbf{A}=\mathbf{D}^{-1}\tilde{\mathbf{W}},
\end{equation}
which is the desired graph Laplacian.
It has been shown that in the limit as $\varepsilon \to 0$ and $m \to \infty$, $\mathbf{A}$ converges to the Laplace-Beltrami operator on $\mathcal{M}$ \cite{Coifman2006}.
Because $\mathbf{A}$ is a Markov matrix, its eigenvalues ($\lambda$) are real-valued and vary between 0 and 1 and its eigenvectors ($\Phi$, stacked column-wise by convention) are real-valued; the trivial eigenvector, $\phi_0=\mathbf{0}$ has eigenvalue $\lambda_0=1$.
We discard $\phi_0$ and sort the remaining eigenvectors $\{\phi_i\}, i=1,...m$ in order of descending eigenvalue.
The first several eigenvectors provide an effective parameterization of the manifold, which call the DMAPS embedding of the data.
Sometimes, to emphasize a spectral gap, we will use $\phi_i\lambda_i^k$ (where $k=1,2,...$) as the DMAPS embedding.

Some of the DMAPS eigenvectors may simply be harmonics that provide no new information about $\mathcal{M}$.
These can be safely discarded either by visual inspection or (more systematically) by using local linear regression to test whether a new eigenvector can be predicted using information from the previous eigenvectors.
When using local linear regression, Dsilva et al.  define a relative leave-one-out cross validation residual for each eigenvector as  \cite{Dsilva2015}:

\begin{equation} \label{eqn:LLR}
R_k=\sqrt{\frac{\sum_{i=1}^n\left(\phi_k(i)-\hat{\alpha}_k(i) + \hat{\beta}_k(i)^T\Phi_{k-1}(i)))\right)^2}{\sum_{i=1}^n\left(\phi_k(i)\right)^2}}
\end{equation}
where $\hat{\alpha}$ and $\hat{\beta}$ are coefficients from regression and $\Phi_{k-1}$ is a matrix containing the first $k-1$ eigenvectors.
If $R_k\ll 1$, then $\phi_k$ can be predicted from the previous eigenvectors and therefore provides only redundant information.

\section{Nonisothermal CSTR Model}\label{app:cstr}
The MPC controller for the CSTR is based on a mechanistic model of a constant density reactor (c.f. \cite{Meadows1996}).
The model equations are:

\begin{equation}\label{eqn:cstr} 
\begin{aligned}
\dot{C}_A &=\frac{q}{V}(C_{A0}-C_A) - kC_A \\
\dot{T_r} &= \frac{q}{V}(T_0-T_r)- \frac{\Delta H}{\rho C_p}kC_A + \frac{UA}{\rho C_P V}(T_c - T_r)\\
C_B &= C_{A0} - C_A
\end{aligned}
\end{equation}
$C_i$ for  $i\in \{A,B\}$ are the species concentrations; $C_{A0}$ is the concentration of species $A$ at the inlet (concentrations of species $B$ is zero at the inlet).
$T_r$ is the reactor temperature, $T_0$ is the inlet temperature, and $T_c$ is the cooling water temperature---$T_c$ is used here as the manipulated variable.
$q$ is the flow rate, $V$ is the reactor volume, $C_p$ is the heat capacity, $\rho$ is the density, $U$ is a heat transfer coefficient, $A$ is the heat transfer area, and $\Delta H$ is the reaction enthalpy (note that $\Delta H$ is a negative number as the reaction is exothermic).
The reaction rate constants, $k$ are given as an Arrhenius expression:

\begin{equation}\label{eqn:arrhenius}
k = k_{0} e^{-\frac{E}{RT_r}}
\end{equation}
$k_{0}$ is the pre-exponential factor, $E$ is the activation energy, $R$ is the gas constant.
Parameter values are given in Table \ref{tab:params}.
To improve numerical performance, we nondimensionalized the concentration variables using $\hat{C}_i=\frac{C_i}{10 \mathrm{~mol~L}^{-1}}$ for $i\in \{A,B\}$ and temperature variables using $\hat{T}_i=\frac{T_i-300 \mathrm{~K}}{100 \mathrm{~K}}$ for $i\in \{r, c\}$.
All of the results shown in Section \ref{sec:results} use the nondimensionalized variables.
\begin{table}
\begin{center}
\caption{Parameter values for the CSTR model.}
\label{tab:params}
\begin{tabular}{| c | c | c| }
\hline
Parameter & Value & Units \\
\hline
 $\frac{E}{R}$ & 6000.0 & \si{\kelvin} \\ 
 $k_{0}$ & $e^{17.5}$ & \si{\per \second} \\ 
 $\frac{\Delta H}{\rho C_p}$ & -16.0 & \si{\kelvin \litre \per \mole} \\
 $\frac{UA}{\rho C_p V}$ & 0.3 & \si{\per \second} \\
 $\frac{q}{V}$ & 1.0 & \si{\per \second}\\
 $C_{A0}$ & 10.0 & \si{\mole \per \litre}\\
 $T_0$ & 300.0 & \si{\kelvin}\\
 \hline
\end{tabular}
\end{center}
\end{table}

\section{MPC Controller}\label{app:mpc}
To build the datasets, $\mathbf{U}^*$, $\mathbf{X}_\alpha^*$, $\mathbf{X}_\beta^*$, and $\mathbf{X}_\gamma^*$, we designed a fully nonlinear MPC controller.
The system was discretized with a sampling time of 0.05 \si{\second}.
The control policy was obtained from the current system state by solving the following optimization problem:

\begin{equation} \label{eqn:mpc_controller}
\begin{aligned}
& \underset{\mathbf{u}} {\text{min}} 
& V(x,k,\mathbf{u}) = \sum_{i=k}^{k+19}{\left(C_{B,i} - r_0 \right)^2 }\\
&\text{s.t.} & u_i \le 2.0,~i=k, k+1, ..., k+19 \\ 
& & u_i \ge -2.0,~i=k, k+1, ..., k+19 \\
& & |u_i - u_{i+1}| \le 0.5,~i=k, k+1, ..., k+19 \\
\end{aligned}
\end{equation}
which indicates a control horizon of 20 time steps (1.0 \si{\second} continuous time), constraints on the maximum and minimum values of $u$, and constraints on the rate of change of $u$.
In Equation \ref{eqn:mpc_controller}, $x=\begin{bmatrix}C_A & T_r \end{bmatrix}$ and $u_k = T_{c,k}$. 
The optimization problem was solved using a sequential quadratic programming algorithm as implemented in the SciPy numerical computing package \cite{Jones,Kraft1988}.

The system was initialized at random states uniformly sampled from $C_A \in [0.1, 0.9]$ and $T \in [0.0, 0.55]$ using random constant references $r_0 \in [0.1, 0.9]$.
The MPC controller was used to evolve the system in time for 20 time steps, and all of the data points from every time step were collected to discover the control manifold and design the explicit MPC controller.

\end{document}